\RequirePackage[l2tabu]{nag}
\documentclass{article}

\usepackage[utf8]{inputenc}
\usepackage[T1]{fontenc}

\usepackage{geometry}
\usepackage{color}
\usepackage{amsmath}
\usepackage{lmodern}
\usepackage[tracking=true,kerning=true,final]{microtype}

\usepackage{hyperref}
\hypersetup{colorlinks=true, linkcolor=blue,  anchorcolor=blue, citecolor=blue, filecolor=blue, menucolor=blue, urlcolor=blue}
\hypersetup{pdfauthor={},pdftitle={}}

\usepackage[caption=false]{subfig} 
\usepackage{color}
\usepackage{graphicx}
\graphicspath{{./}{./figs/}}

\usepackage{amsfonts}
\usepackage{bm}
\usepackage{mathtools}
\usepackage{amsthm}

\newtheorem{theorem}{Theorem}[section]
\newtheorem{proposition}[theorem]{Proposition}
\newtheorem{definition}[theorem]{Definition}

\title{Optimal  criteria and their asymptotic form for data selection in data-driven reduced-order modeling with Gaussian process regression}

\author{Themistoklis P. Sapsis
\thanks{Corresponding author: \href{mailto:sapsis@mit.edu}{sapsis@mit.edu},
Tel: (617) 324-7508, Fax: (617) 253-8689%
} \ and Antoine Blanchard\\
Department of Mechanical Engineering,
\\ Massachusetts Institute of Technology, \\
77 Massachusetts Ave., Cambridge, MA 02139}
\date{\today}

\usepackage{amsmath}
\usepackage{amssymb}

\usepackage[normalem]{ulem}

\begin{document}

 \maketitle\ 
\begin{abstract}
We derive criteria for the selection of datapoints used for data-driven reduced-order modeling and other areas of supervised learning based on Gaussian process regression (GPR). While this is a well-studied area in the fields of active learning and optimal experimental design, most  criteria in the literature are empirical. Here we introduce an optimality condition for the selection of a new input defined as the minimizer of the distance between the approximated output probability density function (pdf) of the reduced-order model and the exact one. Given that the exact pdf is unknown, we define the selection criterion as the supremum over the unit sphere of the native Hilbert space for the GPR. The resulting selection criterion, however, has a form that is difficult to compute. We combine results from GPR theory and asymptotic analysis to derive a computable form of the defined optimality criterion that is valid in the limit of small predictive variance. The derived asymptotic form of the selection criterion leads to convergence of the GPR model that guarantees a balanced distribution of data resources between probable and large-deviation outputs, resulting in an effective way for sampling towards data-driven reduced-order modeling. 

\textbf{Keywords}: Optimal experimental design; Data-driven modeling; Bayesian regression; Optimal sampling; Active learning
  
\end{abstract}

\section{Introduction}
Reduced-order modeling has been a cornerstone of modern computational methods. The effectiveness of the reduced-order models relies both on their design, so they can capture the complexity of the underlying process, and also on the information they rely on. This fundamental information can have the form of  i) governing equations, which typically carry assumptions or simplifications of their own, and ii) data, which is a more reliable but expensive source. The present work involves the development of criteria for the most effective selection of data---or associated experiments to generate this data---in order to perform data-driven reduced-order modeling.

The literature related to data-driven reduced-order modeling is vast and spans a great number of engineering and scientific fields ranging from fluid mechanics \cite{brunton2019a, Karniadakis2021, Ghattas2021,Fernex2020, Ma2021} to structural mechanics \cite{Kerschen2006, Jain2021, Moore2019, Bertalan2019}. In the majority of these works the assumption is plentiful data. While for many applications this is indeed the case, there are several important problems where plentiful data is not available, either because the associated physical or numerical experiments are too expensive or because of the nature of the problem, e.g., extreme events that occur rarely \cite{Sapsis2021}.  

The scientific field that aims to tackle this issue is active learning \cite{Sacks1989, Chaloner95, Gramacy2009}. A critical issue in active learning is the choice of acquisition function, i.e., the criterion used to select which sample to query next in an optimal manner.  Acquisition functions come in various shapes and forms \cite{Chaloner95, Shahriari2016}, but many popular criteria suffer from severe limitations, including high computational cost, intractability in high dimensions, and inability to discriminate between active and idle input variables \cite{sapsis20}. Recently, a new class of active learning criteria was introduced, designed to take into account the effect of the output on the selection of samples \cite{sapsis20, Blanchard_SIAM21}. This new class of output-weighted criteria has led to significantly improved performance in a variety of problems involving uncertainty quantification, Bayesian optimization \cite{blanchard_jcp21, Blanc_turbu}, and decision making \cite{Blanchard2021a, Yang2021}. However, there is no sound theoretical understanding of its favorable properties. 

In this work we rigorously show the optimal properties of output-weighted acquisition functions when it comes to the problem of selecting samples for data-driven reduced-order modeling in the Bayesian context. Specifically, we employ Gaussian process regression (GPR) as the building block to develop Bayesian surrogates or stochastic reduced-order models.  Within the GPR framework, we derive asymptotically optimal acquisition functions, where optimality is defined in the sense of fastest convergence of the probability density function (pdf) describing the quantity of interest. We demonstrate the derived optimal acquisition functions in a mechanical oscillator subjected to high-dimensional stochastic forcing resulting in non-Gaussian heavy tails, as well as the reduced-order modeling of a beam under axial and transverse stochastic loads that result in buckling and bending.     

\section{Problem setup }
Our aim is to build a data-driven reduced-order model that will capture the behavior of an output quantity $y \in \mathbb{R}$ (assumed scalar for simplicity) with respect to an input variable $x \in \mathbb{R}^n$. The input can represent, for example, initial conditions or parameters governing the evolution of a dynamical system, while the output is any quantity of interest that depends on the input variable. We also assume that the input variable has a prescribed probability distribution function,  $p_x({x})$. 

We are given a set of input datapoints $X={\left\{x_i\right\}}_{i=1}^N$ and corresponding output values
\begin{displaymath}
Y=[y(x_1),...,y(x_N)]^{\mathsf{T}},
\end{displaymath} 
and the objective is to identify a criterion that allow us to select the next input point $x_{N+1}\triangleq h$ that will supplement---together with the corresponding output $y(x_{N+1})$---the dataset $\mathcal{D}=\{X,Y\}$. There is a plethora of criteria in the active-learning field (see, e.g., \cite{Chaloner95} for a review), the majority of which are designed to select the next input $h$ that either minimizes the uncertainty of the reduced-order model, i.e., the posterior variance $\sigma^2_y(x|\mathcal{D})$, or maximizes the information content between input and output variables. However, neither of these methods gives appropriate attention on increasing the accuracy of the resulting model in the regions of the input space where it is needed the most: inputs that result in large deviations of the output variable $y$ from its expected value.  

To accomplish this goal we will focus on accelerating the convergence of the model output pdf, $p_y(s|\mathcal{D})$.  In the context of reduced-order modeling this is the appropriate quantity to consider as it encodes the probabilistic information about the large deviations of $y$, while the input variables have been marginalized. Therefore, while criteria based on the posterior variance aim to minimize the error of the model for all or certain input variables, a criterion based on the output pdf explicitly targets regions of the input space that contribute to the output pdf. Moreover, to emphasize regions associated with large deviations of the output variable $y$, we will define convergence in terms of the \textit{logarithm} of $p_y(s|\mathcal{D})$, using the below function to measure the discrepancy between two pdfs $p_1$ and $p_2$:
\begin{equation}
\mathbb{D}(p_1,p_2)=\int_{S_y}\left|\log p_1({s})-\log
p_2({s})\right|\mathrm{d}{s}.
\end{equation}  
In the above, $S_y$ is a finite domain of $y$ (the output variable), over which we aim to build the reduced-order model. Note that this is different from the Kullback--Leibler divergence, which does not give the same emphasis to low probability events (associated with large deviations of $y$). Moreover, it can be easily shown that the above function defines a metric. In what follows, we first provide a quick review of GPR and its basic properties, and subsequently formulate the optimality condition that our sampling criteria will satisfy. 
      
\subsection{Review of Gaussian process regression}
We  use GPR to build a surrogate for the unknown function $y(x):\mathbb{R}^n\rightarrow \mathbb{R}$. The idea is to utilize a Gaussian prior on the function $y$, i.e.,
\begin{equation}
y_0\sim \text{GP}(m(x),k(x,x')),
\end{equation}
with known mean $m(x):\mathbb{R}^n\rightarrow \mathbb{R}$ and covariance function $k(x,x'):\mathbb{R}^{n\times n}\rightarrow \mathbb{R}$ assumed to be positive-definite. We typically set the mean function $m$ to be identically zero. As for the covariance kernel $k$, we often use the Gaussian covariance:
\begin{equation}
k(x,x')=\sigma_k^2\exp\left( -\frac{\left\Vert x-x' \right\Vert^2}{2\lambda^2} \right).
\end{equation} 
 Conditioning the Gaussian process on the available datapoints, we obtain the predictive process
 \begin{align}
y_*\sim \text{GP}(\bar y(x),\bar k(x,x')),
\end{align}
where the predictive mean, $\bar y(x)$ and predictive covariance $\bar k(x,x')$ are given in terms of the datapoints:
\begin{align}\begin{split}
\bar y(x)&=k(x,X)K(X,X)^{-1}Y \\
\bar k(x,x')&=k(x,x')-k(x,X)K(X,X)^{-1}k(X,x'),
\end{split}
\label{gpr00}
\end{align}
where $k(x,X)=[k(x,x_1),...,k(x,x_N)]\in \mathbb{R}^N$ and $K(X,X)=\left\{ k(x_i,x_j) \right\}_{i,j=1}^N \in \mathbb{R}^{N\times N}$.

The native space for GPR schemes is the reproducing kernel Hilbert space corresponding to the kernel $k$, defined as follows:

\begin{definition}
A Hilbert space $\mathcal{H}_k$ of functions $y:\mathbb{R}^n\rightarrow\mathbb{R}$, with inner product $\left\langle \cdot,\cdot \right\rangle_{\mathcal H_k}$, is called the reproducing kernel Hilbert space (RKHS) corresponding to a symmetric, positive-definite kernel $k$ if
\begin{enumerate}
\item 
for all $x\in \mathbb{R}^n$, $k(x,x')$, as a function of its second argument, $x'$, belongs to $\mathcal H_k$; and
\item for all  $x\in \mathbb{R}^n$ and $y \in \mathcal H_k$, $\left\langle y,k(x,\cdot) \right\rangle_{\mathcal{H}_k}=y(x)$.
\end{enumerate}
\end{definition}

\noindent For special choices of kernels $k$, the RKHS can be characterized through its spectrum. Specifically, we have the following theorem for the characterization of the RKHS in the case of Gaussian kernels:

\begin{theorem}[Wendland, 2004, Theorem 10.12, \cite{Wendland2004}]
Let $k(x,x')=\sigma_k^2\exp( -{\left\Vert x-x' \right\Vert^2}/{2\lambda^2}
)$ be the squared-exponential kernel. The corresponding RKHS $\mathcal H_k$ can be written as\begin{align}
\mathcal
H_k=\left\{ y\in L_2(\mathbb{R}^n)\cap C(\mathbb{R}^n):\left\Vert y \right\Vert _{\mathcal H_k}=\frac{1}{c_0}\int|\mathcal{F}[y](\omega)|^2\exp(\lambda ^2\left\Vert \omega \right\Vert^2/2) \, \mathrm{d}\omega < \infty\right\},
\end{align} 
where $c_0$ is a constant that depends on $n$ and $\lambda$, and $\mathcal F$ is the Fourier transform. 
\end{theorem} 
\noindent The above results shows that for any $y \in \mathcal H_k$, the magnitude of its Fourier transform $|\mathcal{F}[y](\omega)|$ decays exponential fast as $|\omega|\rightarrow \infty$ and the speed of decay increases with $\lambda$. Analogous results exist for the case of Mat{\'e}rn kernels \cite{Wendland2004}. 

A fundamental property of GPR schemes is the fact that one can obtain a priori estimates for the accuracy of the surrogate model. In particular, for the case of GPR schemes, we have the following error estimate: 
\begin{proposition}[Stuart \& Teckentrup, 2018, Proposition 3.5, \cite{Stuart2016}]\label{optthm1} Suppose that $\bar y(x)$ and $\bar k(x,x')$ are given by the GPR scheme \eqref{gpr00}. Then \begin{align}
\sup_{\left\Vert y \right\Vert_{\mathcal H_k}=1}|y(x)-\bar y(x)|=\bar k(x,x)^{\frac{1}{2}}\triangleq \bar \sigma (x),
\end{align}
where the supremum occurs when the functions $y(\cdot)$ and $\bar k(\cdot,x)$ are linearly dependent.
\end{proposition}
\noindent Based on this result which involves the accuracy of the surrogate map, we will derive the corresponding results for the pdf of the output variable, $p_y(s|\mathcal{D})$.
\subsection{Optimality condition for data selection}
We formulate an active-sampling criterion that aims directly for the convergence of the output pdf,
 $p_y(s|\mathcal{D})$. Let $h \in \mathbb{R}^n$ be the new candidate input point. As output, we employ the approximation by the surrogate model, $\bar y(h)$. In this way, we have the augmented dataset $\mathcal{D'}=\{[X,h],[Y,\bar y(h)]\}$, which results in the same predictive mean, $\bar y(x)$, and a new predictive covariance $\bar k'(x,x';h)$. Ideally, we would want the new input to be chosen by minimizing the distance
 \begin{equation}
\mathbb{D}(p_y,p_{\bar y'})=\int_{S_y}\left|\log p_y({s})-\log
p_{\bar y'}({s}|\mathcal{D},h)\right| \mathrm{d}{s}.
\end{equation}
However, this is not possible as $p_y$ is a priori unknown. To this end, we will use as a selection criterion the supremum of the above distance over the unit sphere of the functional space $y \in \mathcal{H}_k$ (we fix the norm of the unknown function without loss of generality). Therefore the selection criterion takes the form of minimizing the acquisition function
\begin{equation}
Q(h|\mathcal{D})\triangleq\sup_{ \left\Vert y \right\Vert_{\mathcal{H}_k}=1}\mathbb{D}(p_y,p_{\bar y'})=\sup_{ \left\Vert y \right\Vert_{\mathcal{H}_k}=1}\int_{S_y}\left|\log p_y({s})-\log
p_{\bar y'}({s}|\mathcal{D},h)\right|\mathrm{d}{s}.
\label{criterionq}
\end{equation}
While the criterion is targeting the output pdf, it is not easily computable, especially in high dimensions. The rest of the paper aims to derive a computable version appropriate for high-dimensional input spaces.

\section{Asymptotically optimal criterion for data selection}
Our efforts focus on obtaining a computable version of the criterion \eqref{criterionq}. We plan to achieve this by assuming small variance $\bar \sigma^2(x)=\bar k(x,x)$. We first recall an asymptotic result that connects the error  for a map and the error between the induced pdfs  defined by the corresponding maps. 
\begin{theorem}[Mohamad \& Sapsis, 2018, Theorem 2, \cite{mohamad2018}]\label{approxthm1}
Let $\hat y(x)$ and $y(x):\mathbb{R}^n\rightarrow\mathbb{R}$ be two continuous functions with difference $\Delta  y({x})=\hat y(x)-y(x)$, which is assumed to be small. Let also $p_x({x})$ be the probability density function of the random vector $x \in \mathbb{R}^n$.  The difference between the induced pdfs for $\hat y$ and $y$ has the following asymptotic behavior:
\begin{equation*}
 p_{y}({s})-p_{\hat y}({s})  =-\frac{\mathrm{d}}{\mathrm{d}{s}}\int\displaylimits_{{\hat y}({x})={s}} p_x({x})\Delta  y({x})\,\mathrm{d}{x}+\mathcal{O}(|\Delta  y|^{2}).
\end{equation*}
\end{theorem}
\noindent Building on this result, we have the main theorem that characterizes the supremum in the RKHS between the pdf induced by the surrogate approximation and the pdf induced by the exact map: \begin{theorem}\label{thm_main}
Let $y(x) \in \mathcal{H}_k$ be an arbitrary function and its GPR approximation with kernel $k$   given by ${\bar y(x)}$ with corresponding variance $\bar \sigma^2(x)$. Then the following property holds for small $\bar \sigma$:
\begin{align}\label{res1}
& \sup_{ \left\Vert y \right\Vert_{\mathcal{H}_k}=1}\int_{S_y}\left|\log p_{\bar y}({s})-\log
p_{y}({s})\right|\mathrm{d}{s} =   \int\displaylimits_{\bar y^{-1}(S_y)}\frac{p_x({x})|p'_{ y}({\bar y}({x}))|}
{p_{ y}({\bar y}({x}))^2}\bar \sigma({x})\,\mathrm{d}{x}+\mathcal{O}(\bar \sigma^2).
\end{align}
\end{theorem}
\noindent \textbf{Proof: } Utilizing Theorem \ref{approxthm1}, we have for the case where ${\bar y}$ is close to $y$ (i.e., $\Delta y/ y= (\hat y-y)/y \ll1$)
\begin{align*}\log p_{\bar y}({s})-\log
p_{y}({s})&=\frac{p_{\bar y}({s})-p_{y}({s})}{p_{y}({s})}+\mathcal{O}(|\Delta p_{y}|^2)\\
&=-\frac{\frac{\mathrm{d}}{\mathrm{d}{s}}\int\displaylimits_{{\bar
y}({x})={s}} p_x({x})\Delta  y({x})\,\mathrm{d}{x}}{p_{y}({s})}+\mathcal{O}(|\Delta
 y|^2).
\end{align*}
Expressing the right-hand side as a volume integral with a delta function and using properties of generalized derivatives, we have
\begin{align*}
\frac{\frac{\mathrm{d}}{\mathrm{d}{s}}\int\displaylimits_{{\bar
y}({x})={s}} p_x({x})\Delta  y({x})\,\mathrm{d}{x}}{p_{y}({s})}&=\frac{\frac{\mathrm{d}}{\mathrm{d}{s}}\int\displaylimits_{} p_x({x})\Delta  y({x})\delta
({s-{\bar
y}({x})})\,\mathrm{d}{x}}{p_{y}({s})}\\
&=\int\displaylimits_{}\frac{p_x({x})\Delta
 y({x})\delta '({s-{\bar
y}({x})})}{p_{y}({s})}\,\mathrm{d}{x}\\
&=\int\displaylimits_{} \frac{p_x({x})\Delta
 y({x})p'_{y}({s})\delta ({s-{\bar
y}({x})})}{p^2_{y}({s})}\,\mathrm{d}{x}\\
&=p'_{y}({s})\int\displaylimits_{} \frac{p_x({x})\Delta
 y({x})\delta ({s-{\bar
y}({x})})}{p^2_{y}({s})}\,\mathrm{d}{x}.
\end{align*}
We note that $p_x({x})\delta ({s-{\bar
y}({x})})/p^2_{y}({s}) \ge 0$. By employing Proposition \ref{optthm1} we obtain the tight bound
\begin{align*} \int \frac{p_x({x})\Delta
 y({x})\delta ({s-{\bar
y}({x})})}{p^2_{y}({s})} \, \mathrm{d}{x}\leq   \int\displaylimits_{}
\frac{p_x({x})\bar \sigma({x})\delta ({s-{\bar
y}({x})})}{p^2_{y}({s})}\,\mathrm{d}{x},\end{align*}
where equality holds for the case where the functions $y(\cdot)$ and $\bar k(\cdot,x)$ are linearly dependent \cite{Stuart2016}.

Combining the above, we have
\begin{align*}
\left|\log
p_{\bar y}({s})-\log
p_{y}({s})\right|&\leq  |p'_{y}({s})|\int\displaylimits_{} \frac{p_x({x})\bar \sigma({x})\delta ({s-{\bar
y}({x})})}{p^2_{y}({s})}\,\mathrm{d}{x}+\mathcal{O}(\bar \sigma^2).
\end{align*}
Integrating over ${s}$, we have the final result:
\begin{align*}\sup_{\left\Vert y \right\Vert_{\mathcal{H}_k}=1}\int_{S_y}\left|\log
p_{\bar y}({s})-\log
p_{y}({s})\right| \mathrm{d}s &=  \int_{S_y} |p'_{y}({s})|\int\displaylimits_{} \frac{p_x({x})\bar \sigma({x})\delta
({s-{\bar
y}({x})})}{p^2_{y}({s})}\,\mathrm{d}{x}\,\mathrm{d}{s}+\mathcal{O}(\bar \sigma^2) \\ & =  \int\displaylimits_{\bar y^{-1}(S_y)} \frac{p_x({x})|p'_{y}({{\bar
y}({x})})|}{p^2_{y}({\bar
y}({x}))}\bar \sigma({x})\,\mathrm{d}{x}+\mathcal{O}(\bar \sigma^2).\end{align*}
This completes the proof. $\blacksquare$

Next, we utilize this asymptotic form to reformulate the data selection criterion \eqref{criterionq}. Specifically, we apply the above theorem for the augmented dataset $\mathcal{D}'$ and obtain the asymptotic reformulation of the selection criterion:
\begin{align}
\sup_{ \left\Vert y \right\Vert_{\mathcal{H}_k}=1}\int_{S_y} \left|\log p_y({s})-\log
p_{\bar y'}({s}|\mathcal{D},h)\right|\mathrm{d}{s} & =   \int\displaylimits_{\bar
y^{-1}(S_y)}\frac{p_x({x})|p'_{ y}({\bar y}({x}))|}
{p_{  y}({\bar y}({x}))^2}\bar \sigma({x;h})\,\mathrm{d}{x}+\mathcal{O}(\sigma^2),
\end{align}
where $\bar \sigma({x;h})=\bar k(x,x;h)$ is the predictive variance based on the augmented dataset  $\mathcal{D}'$. We note that the right-hand side involves the pdf $p_y$, which is unknown but can always be approximated by $p_{\bar y}$ with negligible error, given the assumptions of Theorem \ref{thm_main}. This gives us the final asymptotic approximation:
\begin{align}
Q(h|\mathcal{D})=\sup_{ \left\Vert y \right\Vert_{\mathcal{H}_k}=1}\int_{S_y}\left|\log p_y({s})-\log
p_{\bar y'}({s}|\mathcal{D},h)\right|\mathrm{d}{s} & \simeq
\int\displaylimits_{\bar y^{-1}(S_y)}\frac{p_x({x})|p'_{\bar y}({\bar y}({x}))|}
{p_{\bar  y}({\bar y}({x}))^2}\bar \sigma({x;h}) \,\mathrm{d}x.
\end{align}
It is worth emphasizing the term in the denominator, which promotes sampling of regions associated with low probability, i.e., large deviations of the output of the reduced-order model. 

We can simplify the right-hand side further by using the Cauchy--Schwarz inequality (with weight $\frac{p_x({x})}{p_{\bar  y}({\bar
y}({x}))}$) to obtain the less conservative upper bound:
\begin{equation}\label{res2}
Q(h|\mathcal{D})\simeq\int\displaylimits_{\bar y^{-1}(S_y)} \frac{p_x({x})|p'_{\bar y}({{\bar
y}({x})})|}{p^2_{\bar y}({\bar y}({x}))}\bar \sigma({x;h})\,\mathrm{d}{x}\le c  \left[\int\displaylimits_{\bar y^{-1}(S_y)}\frac{p_x({x})}{p_{\bar  y}({\bar y}({x}))}\bar \sigma^{2}({x};h)\, \mathrm{d}x\right]^\frac{1}{2},
\end{equation}
where 
\begin{equation}
c=\left[\int\displaylimits_{\bar y^{-1}(S_y)} \frac{p_x({x})p^{\prime 2}_{\bar y}({{\bar
y}({x})})}{p^3_{\bar y}({\bar y}({x}))}\,\mathrm{d}{x}\right]^\frac{1}{2}
\end{equation} is a constant that depends only $p_x(x)$, $p_{\bar y}(y)$, and $\bar
y(x)$, i.e., not on $h$.  This form, also referred as output-weighted (or likelihood-weighted) criterion, is appropriate for computations involving even high-dimensional input spaces, since it allows for the analytical computation of $\sigma^{2}({x};h)$ in terms of simpler integrals \cite{Blanchard_SIAM21}. It has been studied numerically in recent papers \cite{sapsis20, Blanchard_SIAM21,blanchard_jcp21}, showing significantly favorable convergence properties compared with existing active-learning criteria.

\subsection{The case of extreme-event quantiles}
For a wide range of applications, the focus is on characterizing the probability of exceeding a certain level, $P_y(s_*)=P[y\leq s_*]$, rather than characterizing the full pdf $p_y(s)$  (the case of non-exceeding probability can be addressed in a similar fashion).  For this type of problems the appropriate selection criterion is minimizing the acquisition function\begin{equation}
R(h|\mathcal{D},s_*)\triangleq\sup_{ \left\Vert y \right\Vert_{\mathcal{H}_k}=1}| P_{\bar y'}({s_*}|\mathcal{D},h)-
P_{y}({s_*})|.
\end{equation}For this case, we have the following asymptotic result for small $\sigma$:
\begin{theorem}
Let $y(x) \in \mathcal{H}_k$ and its GPR approximation with kernel $k$  
given by ${\bar y(x)}$ with variance $\bar \sigma^2(x)$. Then, for any given $s_*$, the
following property holds for small $\bar \sigma$:
\begin{align*}
& \sup_{ \left\Vert y \right\Vert_{\mathcal{H}_k}=1}| P_{\bar y}({s_*})-
P_{y}({s_*})| =   \int\displaylimits_{\bar
{
y}({x})={s_*}}\bar \sigma({x})p_x({x})\,\mathrm{d}x+\mathcal{O}(\bar \sigma^2).
\end{align*}
\end{theorem}
\noindent \textbf{Proof: } The starting point is Theorem \ref{approxthm1}. We integrate from $-\infty$ to $s_*$ to obtain the cumulative distribution function on the left-hand side, when ${\bar y}$ is close
to $y$ (i.e., $\Delta y/ y= (\hat y-y)/y \ll1$):
\begin{align*}P_{\bar y}({s})-P
_{y}({s})=-\int\displaylimits_{{\bar y}({x})={s}} p_x({x})\Delta  y({x})\,\mathrm{d}{x}+\mathcal{O}(|\Delta
 y|^2).
\end{align*}
We  then employ Proposition \ref{optthm1} and bound the difference $\Delta  y$ on the right-hand side. This completes the proof.  $\blacksquare$

Applying the above theorem to the augmented dataset $\mathcal D'$, we obtain the asymptotic form for the optimal selection criterion involving quantiles:
\begin{align}
& R(h|\mathcal{D},s_*)\triangleq\sup_{ \left\Vert y \right\Vert_{\mathcal{H}_k}=1}| P_{\bar y'}({s_*}|\mathcal{D},h)-
P_{y}({s_*})|=   \int\displaylimits_{\bar
{
y}({x})={s_*}}\bar \sigma({x;h})p_x({x})\,\mathrm{d}x+\mathcal{O}(\bar \sigma^2).
\end{align}

We observe that for this case, the optimal data selection criterion takes a different form, focusing primarily on reducing the error around the contour $\bar y(x)=s_*$. This is not surprising given that the extreme-event quantile is essentially equivalent to a classification problem, so what is most important is to have low error around the a priori unknown contour of interest. 

\subsection{Convergence of spatial approximation error}
Given an active-learning criterion, it is important to characterize the convergence of the spatial error of the GPR approximation. This can be obtained with the help of Proposition \ref{optthm1} and standard results of measurable functions. Here we give a proof for the convergence properties using the sampling criterion appearing on the right-hand side in \eqref{res2}. Specifically, we have the following result:
  
\begin{theorem}
Suppose that ${\bar y_N}$ and $\bar \sigma_N$ are given by a GPR with kernel $k$, using $N$ samples. Moreover, assume that sampling is performed using the optimal criterion so that \begin{displaymath}
\lim_{N\rightarrow\infty} \int \frac{p_x({x})}
{p_{\bar y_N}({\bar y_N}({x}))}{\bar \sigma_{{N}}^2}({x})\,\mathrm{d}{x}=0.
\end{displaymath} 
Then we have convergence in measure, i.e.,
\begin{equation*}
\lim_{N\rightarrow \infty}P \left[ x:\sup_{ \left\Vert y \right\Vert_{\mathcal{H}_k}=1}({\bar y_N}({x})-{y}({x}))^2 \le q \frac{p_{\bar y_N} ({\bar y_N}( x))}{p_x({x})} \right]=1,
\end{equation*}
for every $q>0$.
 \end{theorem}
 \textbf{Proof:} From standard results of measurable functions \cite{vulikh}, for every sequence of functions $\phi_k(x)\in L^p$ ($1\le p \le\infty$) for which $\Vert \phi_k(x) \Vert_p \rightarrow  0$, we have convergence in measure:
 \begin{align} 
\lim_{k\rightarrow \infty} P[x:|\phi_k(x)|\ge q]=0,
 \end{align} 
 for every $q>0$. Applying this result to the sequence $p_x({x})\bar \sigma_{{N}}^2({x})/ 
p_{\bar y_N}({\bar y_N}({x}))$, we obtain
\begin{align*}\lim_{N\rightarrow \infty} P \left[x:\frac{p_x({x})}
{p_{\bar y_N}({\bar y_N}({x}))}{\bar \sigma_{{N}}^2}({x})\le q \right]=1.
\end{align*}
We then utilize Proposition \ref{optthm1}, which immediately leads to the desired result. $\blacksquare$

This result provides a description of the spatial convergence properties for the approximation error of the reduced-order model. Specifically, it shows that the convergence is accelerated in regions of the input space that are i) most probable according to the pdf $p_x(x)$, and ii) associated with small probability of the output pdf $p_{\bar y_N} ({\bar y_N}( x))$, i.e., large deviations of the output $y$. This type of convergence guarantees a balanced distribution of resources between the most probable inputs and those that result in large deviations for the output, resulting in an effective way of sampling towards data-driven reduced-order modeling.

\section{Numerical illustration}

We demonstrate the optimal sampling criteria in a mechanical oscillator subject to stochastic forcing and in the reduced-order modeling of a beam under axial and transverse stochastic loads.  We consider the optimal sampling criteria appearing on each side of the inequality in \eqref{res2}.  The left-hand side will be referred to as the ``B'' criterion, and the right-hand side as ``IVR-LW'' as it is strictly equivalent to the eponym criterion introduced in \cite{sapsis20, Blanchard_SIAM21}.  We also consider two criteria commonly used in the literature which do not account for the importance of the output relative to the input; namely, uncertainty sampling (US) and input-weighted integrated variance reduction (IVR-IW), whose definitions can be found in \cite{sapsis20, Blanchard_SIAM21}.  For each example below, we run 100 Bayesian experiments, each differing in the choice of the $n+1$ points making up the initial dataset.  Observations are assumed to be corrupted by Gaussian noise with zero mean and unknown (i.e., to be learned) variance $\sigma_n^2$.  Performance at each iteration is evaluated using the median of $\mathbb{D}(p_y,p_{\bar{y}_N})$ across the 100 randomized experiments. Note that we employ this particular quantity to measure performance so we can better emphasize the accuracy in the tails of the resulting probability distributions.

\subsection{Forced nonlinear oscillator exhibiting extreme events}

We begin with the stochastic oscillator of Mohamad and Sapsis \cite{mohamad2018},
\begin{equation}
\ddot{u} + \delta \dot{u} + F(u) = \xi(t), \quad t \in [0,T],
\label{eq:43}
\end{equation}
where $u(t) \in \mathbb{R}$ is the state variable, $F$ a nonlinear restoring force, and $\xi(t)$ a stationary stochastic process which we parametrize using a Karhunen--Lo{\`e}ve expansion with $n$ modes:\begin{displaymath}
\xi(t)\simeq x \Phi(t),\ 
\end{displaymath}  where $\{\Lambda, \Phi(t)\}$ contains the first $n$ eigenpairs of the correlation matrix, and $x \in \mathbb{R}^n$ is a vector of random coefficients having mean zero and diagonal covariance matrix $\Lambda$.  The quantity of interest is taken to be the mean value of $u(t)$ over the interval $[0,T]$.  (Further details about system parameters can be found in \cite{mohamad2018, Blanchard_SIAM21}).

\begin{figure}[!ht]
\centering 
\includegraphics[width=5.9in, clip=true, trim=15 10 10 10]{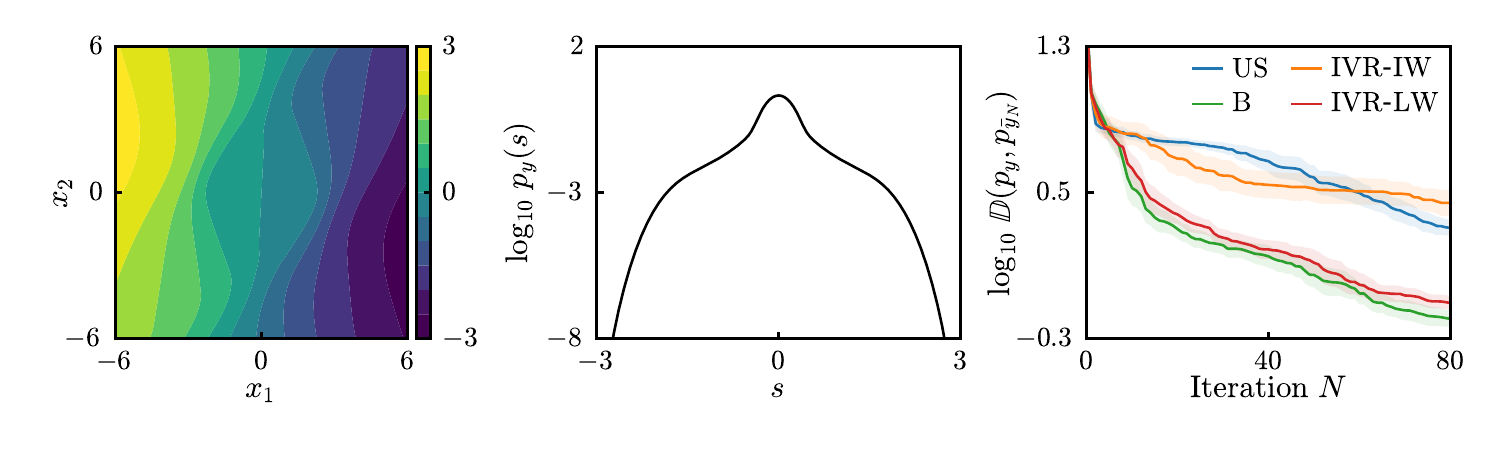}
\caption{For the stochastic oscillator \eqref{eq:43} with $n=2$, contour plot of the output $y$ (left) and its pdf (center), and performance of several sampling criteria for $\sigma_n^2 = 10^{-3}$ (right). The error bands indicate one half of the median absolute deviation.}
\label{fig:1}
\end{figure}

We consider the case $n=2$ as it allows visual comparison of the decisions made by the sampling criteria as more points are being acquired.  Even with $n=2$, the output pdf has heavy tails (see figure \ref{fig:1}), a consequence of the strong nonlinearity in \eqref{eq:43}.   For $\sigma_n^2 = 10^{-3}$, figure \ref{fig:1} shows that the derived optimal criteria accelerate convergence of the output pdf quite dramatically.  The error for B remains close to, but always slightly below, that for IVR-LW, which is consistent with the mathematical derivation laid out in the previous section.  

\begin{figure}[!t]
\centering 
\subfloat[][US]{\includegraphics[width=2.2in]{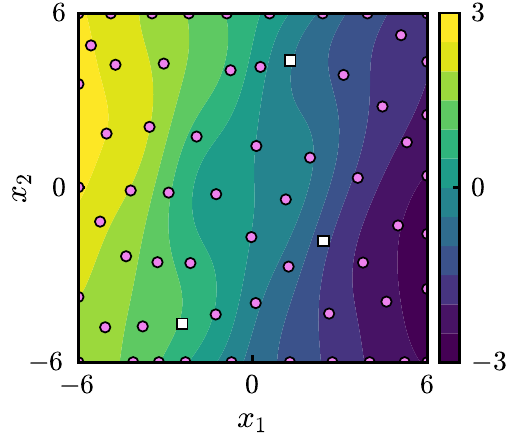}} \qquad
\subfloat[][IVR-IW]{\includegraphics[width=2.2in]{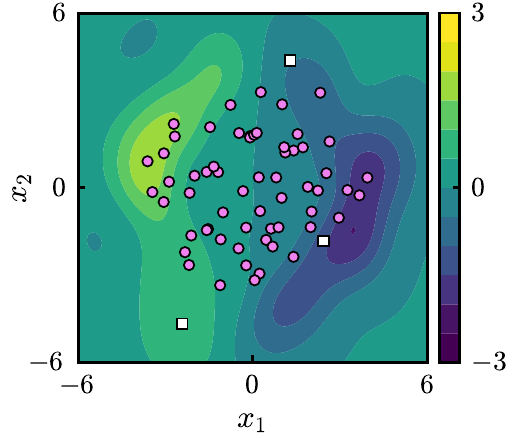}}

\subfloat[][B]{\includegraphics[width=2.2in]{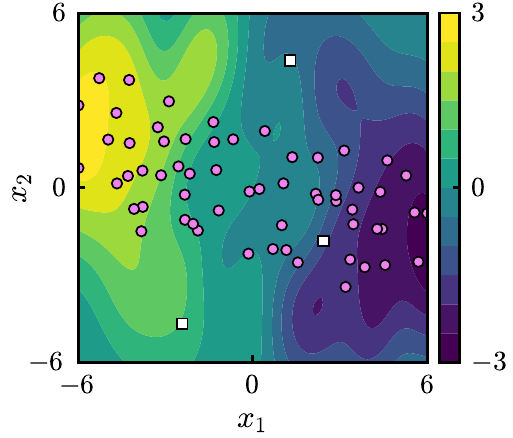}} \qquad
\subfloat[][IVR-LW]{\includegraphics[width=2.2in]{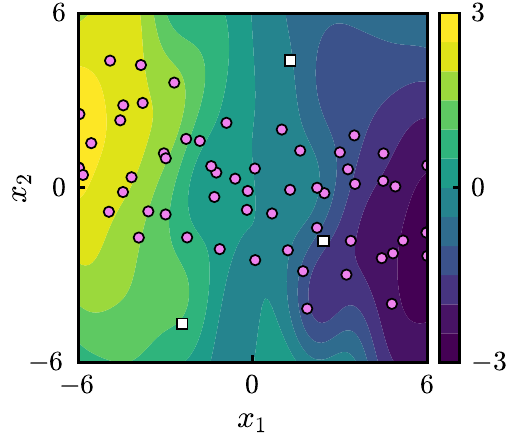}}

\caption{For the stochastic oscillator \eqref{eq:43} with $n=2$ and $\sigma_n^2=0$, progression of the sampling algorithm for several criteria after 60 iterations; the contours denote the posterior mean of the GPR model, the open squares the initial dataset, and the filled circles the optimized samples.}
\label{fig:2}
\end{figure}

To explain the success of the proposed optimal criteria, we investigate the decisions made by US, B, and IVR-LW in the case where observations are noiseless and $\sigma_n^2$ is set to zero in the GPR model (i.e., it is not learned from data).  Consistent with \cite{sapsis20, Blanchard_SIAM21}, figure \ref{fig:2} shows that US attempts to reduce uncertainty somewhat evenly across the space as it has no mechanism to discriminate between relevant and irrelevant regions.  With IVR-IW, the algorithm does not explore beyond the center region where $p_x$ is large, and consequently the interesting regions are not visited.  On the other hand, both B and IVR-LW decide to focus on a diagonal band, with the former being even more surgical and localized than the latter. Specifically, IVR-LW focuses on input regions with important probability as well as those input regions associated with large outputs.  In this way the resulting surrogates predict the output statistics much better than with US.

\subsection{Buckling of a beam under stochastic axial and transverse excitation}

Next, we consider the case of a beam of length $l$ subject to both axial and transverse stochastic loads. The linearized equation of motion takes the form
\begin{align}
        \frac{\partial^2w}{\partial t^2}+2\zeta\omega_0\frac{\partial w}{\partial t}+\omega_0^2\frac{\partial ^4w}{\partial x^4}+P(t)\frac{\partial ^2w}{\partial x^2}=R(x,t),
\end{align}
where $P(t)$ is random axial load and $R(x,t)$ is a random distributed transverse load \cite{nayfeh_mook}. We assume that both functions have zero-mean Gaussian statistics and prescribed spectra. We also assume that the beam has pin boundary conditions at both ends. In this case and assuming weak damping, the solution can be expressed as
\begin{align}
        w(x,t)=\sum_{j=1}^{\infty}f_j(t)\sin(j\pi x/l),
\end{align}
which results in a set of modal equations:
\begin{align}
        \ddot f_j+2\zeta\omega_0 \dot f + \omega_j^2 [1-P(t)/c_j] f=R_j(t),  \qquad j=1,2,...
\end{align}
where
\begin{align}
\omega_j^2=\omega_0^2 (j\pi/l)^4,  \quad c_j=(j\pi/l)^4, \quad \text{and} \quad R_j(t)=\frac{2}{l}\int_0^l \sin (j\pi x/l)R(x,t)\,\mathrm{d}x.
\end{align}
We define as quantity of interest the maximum absolute displacement at $x=l/4$ over a prescribed time interval $[0,T]$:
\begin{align}
y=\max_{t\in [0,T]} \left|w(l/4,t)\right|= \max_{t\in [0,T]} \left|\sum_{j=1}^{J}\sin(j\pi/4)f_j(t)\right|.
\end{align}
We use $\sigma_\xi^2 \exp[-t^2/(2\ell_{\xi}^2)]$ for the correlation function of $P(t)$ and each $R_j(t)$, which are expanded using Karhunen--Lo{\`e}ve expansions with $n_\textit{KL}$ modes.  The search space, therefore, has dimension $n= J(n_\textit{KL}+1)$.  We use parameters $\sigma_\xi=20$, $\ell_{\xi}=0.1$, and $T=5$.

\begin{figure}[!ht]
\centering 
\includegraphics[width=5.9in, clip=true, trim=12 10 10 10]{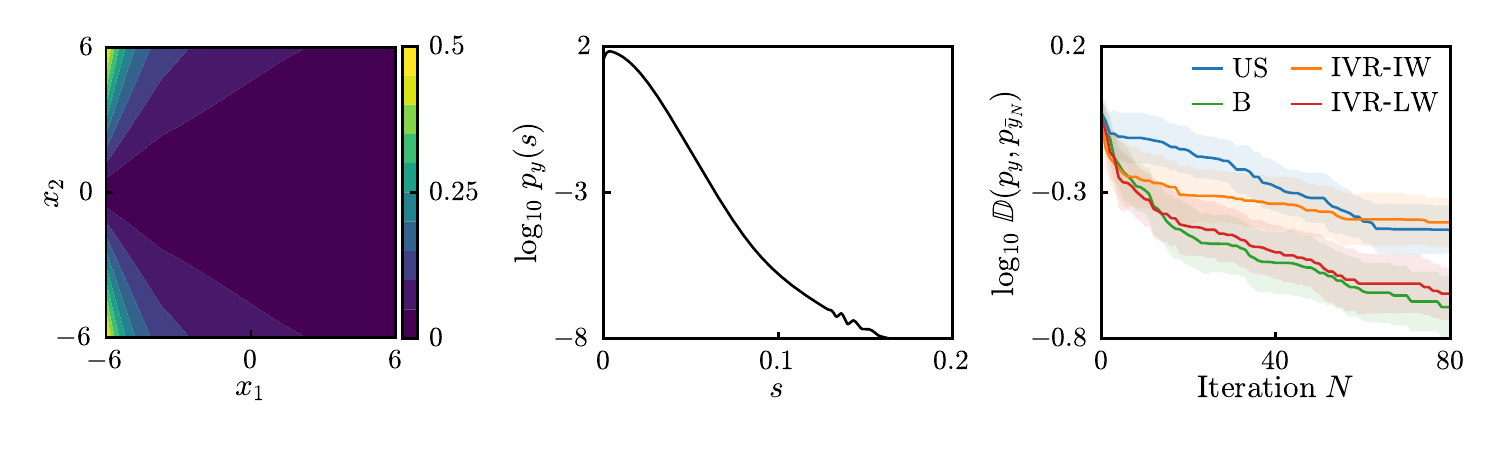}
\caption{For the beam under random load with $n=2$ ($J=1$, $n_\textit{KL}=1$), contour plot of the output $y$ (left) and its pdf (center), and performance of several sampling criteria for $\sigma_n^2 = 10^{-3}$ (right). The error bands indicate one half of the median absolute deviation.}
\label{fig:3}
\end{figure}

For simplicity, we consider a modal truncation of $J=1$ and $n_\textit{KL}=1$, leading to a two-dimensional search space.  For the parameters considered, figure \ref{fig:3} shows that the pdf of the output has a heavy right tail, to which the B and IVR-LW criteria converge more quickly than US and IVR-IW.  The reason is that the extreme displacement values are found in low-probability areas of the search space (i.e., small $p_x$), which US and IVR-IW have no mechanism to discover.

\section{Conclusions}

We have derived optimal acquisition functions for active-learning schemes utilized for reduced-order modeling based on Gaussian process regression. The key feature of these optimal acquisition functions is a mechanism that targets the output pdf of the surrogate model.  The derivation begins by selecting each sample so that the distance between the exact output pdf and the approximated output pdf obtained from the reduced-order model is minimized. Given that the exact pdf is a priori unknown, a supremum of this distance over the native Hilbert space for the Gaussian process regression is considered. The resulting selection criterion, although optimal, it is generally computationally intractable. We show that this difficulty can be overcome by deriving successive asymptotic upper bounds, resulting in a sampling criterion which can be evaluated analytically along with its gradients.  In addition to the data-selection criterion, we derive the corresponding bound for the approximation error of the resulted reduced-order model. This result shows that our approach enables the reduced-order model to find the optimal trade-off between most likely and most interesting (i.e. large-deviation) output for the process at hand, thereby greatly enhancing the effectiveness of the algorithm. Numerical results confirm the derived analytical findings. 

Note that in this work and for the sake of simplicity we have chosen to represent the reduced-order model in the form of a function from the parameter space to the output space. The present framework can be adapted to represent the reduced-order model in the form of a low-dimensional dynamical system, focusing e.g. on capturing a specific mechanism of the full dynamics. We leave this topic as a possible direction for future work.

\subsubsection*{Acknowledgments} 
The authors acknowledge support from the Air Force Office of Scientific Research (MURI\ Grant No. FA9550-21-1-0058),  the Defense Advanced Research Projects Agency (Grant No. HR00112110002), and the Office of Naval Research (Grant No. N00014-21-1-2357).

 \bibliographystyle{abbrv}
\bibliography{library}

\end{document}